\documentclass{article}
\usepackage{amsthm,amsfonts,amsmath,amssymb}
\newcommand{\eq}{\begin{equation}}
\newcommand{\en}{\end{equation}}

\newcommand{\giv}{\,|\,}

\newcommand{\prob}{\mathbb P}
\newcommand{\ex}{\mathbb E}

\newcommand{\ed}{ \stackrel{d}{=}}

\newtheorem{theorem}{\large Theorem}
\newtheorem{proposition}[theorem] {\large Proposition}

\newtheorem{corollary}[theorem]{\large Corollary}
\newtheorem{lemma}[theorem]{\large Lemma}

\begin{document}
\title{Counting the Chain Records: The Product Case}
\author{Alexander V. Gnedin\thanks{Postal address:
 Department of Mathematics, Utrecht University,
 Postbus 80010, 3508 TA Utrecht, The Netherlands. E-mail address: gnedin@math.uu.nl}}
\date{}
\maketitle

\begin{abstract}
\noindent
Chain records is a new type of  multidimensional record.
We discuss how often the chain records are broken when
the background sampling is from the unit cube with uniform 
distribution (or, more generally, from an arbitrary continuous 
product distribution).
\end{abstract}

\section{Introduction}

Consider independent marks  $X_1,X_2,\ldots$  sampled from the
uniform distribution in
$Q_d=[0,1]^d$.
We  define a mark $X_n$ to be a {\it chain record}
if $X_n$ beats the last chain record in $X_1,\ldots,X_{n-1}$.
More precisely, {\it record values} and {\it record  indices} are introduced recursively, 
by setting $T_1=1,\, R_1=X_1$ and
 $$T_{k}=\min\{n>T_{k-1}:~ X_n\prec R_{k-1}\}\,,~~R_k=X_{T_k}\,,~~~~~~k>1\,.$$
Here, $\prec$ denotes the standard strict partial order on ${\mathbb R}^d$ defined in terms of component-wise orders
by
$$
x=(x^{(1)},\ldots, x^{(d)}) \prec y=(y^{(1)},\ldots,y^{(d)})~~
{\rm iff~~}x\not=y~~{\rm and}~~   x^{(i)}\leq y^{(i)}\,~~{\rm for}~~i=1,\ldots,d\,.
$$
It is easy to see that, in any dimension $d$, 
the terms of $(T_k)$ are indeed well defined for all $k$, that is the chain records 
occur  infinitely many times.

\par Although the definition is an obvious restatement of the classical definition
of lower record, this notion of multidimensional record
has not been explored so far. The chain records interpolate between 
two other types of multidimensional records which have been studied in some depth 
\cite{Barn, Bar, DZ, GnDep, GnCone, JMVA, GnNorm,   GR, GRMany, GRScat, HH}.
We  say that a {\it strong record} occurs at index $n$ if either $n=1$, or $n>1$ and 
$$
X_n  \prec X_j ~~{\rm for ~~}j=1,\ldots,n-1.
$$
In the terminology of partially ordered sets, a strong record $X_n$ is 
the {\it least} element in the point set $\{X_1,\ldots,X_n\}$.
Since repetitions in each component have probability zero, $X_n$
is a strong record if and only if there are $d$ {\it marginal} 
strict lower records at index $n$ simultaneously.
We say  that a {\it weak record} occurs at index $n$ if either $n=1$, or $n>1$ and
$$
X_j  \not\prec X_n ~~{\rm for ~~}j=1,\ldots,n-1.
$$
A weak record $X_n$ is a {\it minimal} element in the set $\{X_1,\ldots,X_n\}$. 
Obviously, each strong record is a chain record.
Also, each chain record is a weak record, 
as follows easily by induction from transitivilty of the relation $\prec$.
To illustrate,
for the two-dimensional configuration of points
in Figure 1
the weak records occur at times $1,\,2,\,3,\,5,\,6,\,7,\,8,\,9$, a sole strong record occurs at 1, the marginal records occur at 
1,\,2,\,3,\,5,\,6, and the chain records occur at indices $1,\,5,\,8$.  
\begin{figure}
\unitlength 1mm 
\linethickness{0.4pt}
\ifx\plotpoint\undefined\newsavebox{\plotpoint}\fi 
\begin{picture}(127.54,63.11)(0,0)
\put(67.79,5.48){\framebox(59.75,57.63)[cc]{}}
\put(95.99,49.32){\circle*{1.27}}
\put(78.49,30.23){\circle*{1.27}}
\put(72.83,16.09){\circle*{1.27}}
\put(105.71,19.44){\circle*{1.27}}
\put(87.33,27.05){\circle*{1.27}}
\put(109.78,26.16){\circle*{1.27}}
\put(81.67,11.84){\circle*{1.27}}
\put(117.56,42.25){\circle*{1.27}}
\put(74.6,58.69){\circle*{1.27}}
\put(85.92,51.98){\circle*{1.27}}
\put(121.09,55.68){\circle*{1.27}}
\put(95.81,49.14){\framebox(.18,.35)[]{}}
\put(97.23,51.09){\makebox(0,0)[cc]{1}}
\put(118.97,44.72){\makebox(0,0)[cc]{2}}
\put(107.13,21.92){\makebox(0,0)[cc]{3}}
\put(122.15,58.16){\makebox(0,0)[cc]{4}}
\put(79.37,32.7){\makebox(0,0)[cc]{5}}
\put(76.72,60.81){\makebox(0,0)[cc]{6}}
\put(88.56,29.17){\makebox(0,0)[cc]{7}}
\put(73.54,18.74){\makebox(0,0)[cc]{8}}
\put(84.5,13.43){\makebox(0,0)[cc]{9}}
\put(87.86,53.21){\makebox(0,0)[cc]{10}}
\put(111.37,28.11){\makebox(0,0)[cc]{11}}
\put(67.8,5.57){\framebox(28.06,43.83)[cc]{}}
\put(67.8,5.57){\framebox(10.72,24.6)[cc]{}}
\put(67.8,5.57){\framebox(5.05,10.62)[cc]{}}
\end{picture}
\caption{Chain records in the square} 
\end{figure}
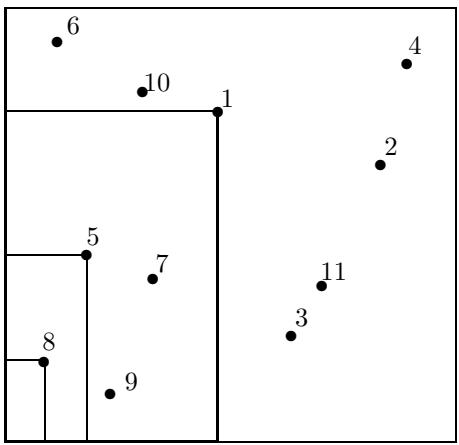
Notably, the chain records are more sensible to arrangement of 
marks in sequence: a permutation of $X_1,\ldots,X_{n-1}$ may destroy or create a chain record at index $n$.  

\par Denote
$\underline{N}_n,\overline{N}_n$ and $N_n$, respectively,  the counts  of strong, 
weak and chain records among the first 
$n$ marks. Thus
$$
\underline{N}_n \leq N_n\leq \overline{N}_n\,.
$$
To underscore concretely the extent of  compromise between weak and strong records,
we need some estimates 
of  how often the  records of different kinds may occur.

\par Recall that in the case $d=1$ the occurences of records are 
independent, with probability $1/n$ for index $n$; this basic fact (known as the Dwass-R{\'e}nyi lemma \cite{Nevzorov,
Resnick})
 implies that
the number of classical records is asymptotically Gaussian with both mean and variance
about $\log n$. This translates easily to the marginal records in $d$ dimensions, since the
marginal rankings are independent. The latter kind of independence 
is characteristic for sampling from 
product distributions in ${\mathbb R}^d$ with continuous marginals, 
hence the instance of $Q_d$ with uniform distribution 
covers the general product case.

\par Properties of the strong-record counts for sampling from $Q_d$ are also rather simple. 
By independence of marginal rankings 
we have a representation $\underline{N}_n=I_1+\ldots+ I_n$ with independent 
Bernoulli indicators and $\underline{p}_n:=
\prob(I_n=1)=n^{-d}$. Thus
$${\mathbb E}\,\underline{N}_n=\sum_{j=1}^n {1\over j^d}\,.$$
Since for $d>1$ the series $\sum\underline{p}_n$ converges, the total number of strong records in the infinite 
sequence of marks is almost surely finite. 

\par Counting the weak records 
is a more delicate matter since their occurences are not independent.
However, we may exploit a correspondence between
weak records in $Q_d$ and the minimal elements in $Q_{d+1}$ (depending on the  context these points are also called 
 Pareto, admissible, efficient, etc.).
The correspondence is established by arranging the marks in $d+1$ dimensions by increase in one fixed component.
By induction in $d$ one can show that 
$${\mathbb E}\,\overline{N}_n=\sum_{1\leq j_1\leq\ldots\leq j_{d}\leq n} 
{1\over j_1\cdots j_d}\sim {1\over d!}\,{(\log n)^d}\,,$$
see \cite{Barn}.
From further known results (see \cite{Bar} and references therein) 
follows that 
the variance ${\rm Var\,}[{\overline N}_n]$ is of the same order $(\log n)^d$, and that
${\overline N}_n$ is asymptotically Gaussian.

\par Thus the strong records are much more rare and the weak records are much more frequent than the
classical records.
In this note we show that, as far as the frequency is concerned,
the chain records in any dimension $d$ are more  in line with the classical records:

\begin{proposition}\label{P1}
For sampling from $Q_d$ with uniform distribution the number of chain records
$N_n$ is approximately Gaussian with moments
$${\mathbb E}\,[N_n]\sim d^{-1}\log n\,,~~{\rm Var}\, [N_n]\sim d^{-2}\log n\,.$$
\end{proposition}
\noindent
The CLT will be proved in Section 3. Above that, we will derive exact and asymptotic formulas
for the probability of a chain record and discuss some scaling limits.

\par The chain records comprise a `greedy' chain in $\prec$, meaning  that  
a mark is joined each time the chain constraint is not violated.
More efficient nonanticipating algorithms for constructing
long chains were designed in 
\cite{BarGn},
 and the length
of the longest possible chain on $n$ random marks was estimated in 
\cite{BW}.  From yet another perspective, the sequence of chain records corresponds to a particular path
in a random data structure called {\it quad-tree} \cite{Devroye, Flajolet}.

\section{The  heights at records}

For $x\in Q_d$ the quadrant $L_x:=\{y\in Q_d:\, y\prec x\}$ is the  lower section of the partial order at $x$.
 The {\it height} $h(x)$ is 
the product of coordinates,
which in the case of uniform distribution under focus is equal to the value of the multidimensional 
distribution function at $x$, i.e. 
the measure of $L_x$.
The height is a key quantity to look at, because the heights at chain records determine the sojourns.
Let $H_k=h(R_k)$. 

\begin{lemma}\label{soj0} Given $(H_k)$ the sojourns $T_{k+1}-T_{k}$ are conditionally independent, geometric
with parameters $H_k$, $k=1,2,\ldots$
\end{lemma}
\begin{proof} A new chain record $R_{k+1}$ occurs as
soon as $L_{R_k}$ is hit by some mark.
\end{proof}
\noindent
The lemma has the following elementary but important consequence.

\begin{corollary}\label{soj} Given $(H_k)$, the conditional law for occurencies of the chain records  for any $d$
is the same as in the classical case $d=1$.
\end{corollary}

\par The heights at records undergo  a multiplicative renewal process, sometimes called
{\it stick-breaking}.
 Let   $W, W_1, W_2, W_3,\ldots$ be i.i.d. copies of 
$H_1=h(X_1)$.

\begin{lemma} 
The heights $(H_k)$ have the same law as the sequence of products
$(W_1\cdots W_k,~\,k=1,2,\ldots)$.
\end{lemma}
\begin{proof}Each lower section $L_x$, viewed as a partially ordered probability space with
normalised Lebesgue measure is isomorphic to $Q_d$ (via a coordinate-wise scale transformation).
Hence all  ratios $H_{k+1}/H_k$ are i.i.d., with the same law as $H_1$.
\end{proof}

\noindent
Explicitly, the density of $W$ is 
\eq\label{H}
\prob(W\in {\rm d}s)={(\log s)^{d-1}\over (d-1)!}\,{\rm d}s\,,~~s\in [0,1],
\en
and its Mellin transform is 
$$g(\lambda):= \ex\,[W^\lambda]=(\lambda+1)^{-d} \,,$$
as follows by noting that $H_1$ is the product of $d$ independent uniform variables.

\par The distinction with the classical $d=1$ case is seen already at this early stage of our discussion.
In the classical case $H_1$ has uniform distribution, hence  the stick-breaking sequence 
$(W_1\cdots W_k,~\,k=1,2,\ldots)$
is the sequence of points of a self-similar (i.e. invariant under homotheties) 
Poisson process with intensity ${\rm d}s/s, ~s\in [0,1]$.
For $d>1$ the point process $(H_k)$ is neither Poisson nor self-similar, 
which is a major source of difficulties leading, e.g., to dependencies in the
occurences of chain records at distinct $n$.

\section{Proving the CLT}

Corollary \ref{soj} suggests to focus on properties of a univariate sequence of random variables 
modified by conditioning and then mixing over some 
given distribution for its subsequence of record values.

\par Let $(U_j)$ be a sequence of $[0,1]$ uniform points, independent of 
$(H_k)$. We shall produce a transformed sequence $(U_j)\giv(H_k)$ by  replacing 
some of the terms in $(U_j)$  by the 
$H_k$'s.
Replace $U_1$ by $H_1$. Do not alter $U_2,U_3,\ldots$ as long as they do not hit $[0,H_1[\,$;
then replace the first uniform point hitting the interval $[0,H_1[\,$ by $H_2$. 
Inductively, as  $H_1,\ldots,H_k$ got inserted, 
keep on screening uniforms until first hitting $[0,H_k[\,$,
then insert $H_{k+1}$ in place of the uniform point that caused the hit, and so on. 
Eventually all $H_k$'s will enter the resulting sequence.
It is easy to see that given $(H_k)$ the distribution of $(U_j)\giv(H_k)$ is the same as 
the conditional distribution of $(U_j)$ given the subsequence of record values $(H_k)$.
In the classical case, $(H_k)$ is the stick-breaking sequence with uniform factors, and  we have
$(U_j)\giv(H_k)\ed (U_j)$, so the insertion  does not alter the law of the sequence.

\par By Corollary \ref{soj}, $N_n$ can be identified with the number of points among $U_1,\ldots,U_n$ 
that 
get replaced by some $H_k$'s.
\par 
There is yet another related interpretation
in terms of {\it partially exchangeable partitions}, as introduced in \cite{PTRF}.
The unit interval $]0,1[$ is divided by $(H_k)$
in infinitely many disjoint
subintervals $[H_1,H_0[\,,\,[H_2,H_1[\,,\ldots$ (where $H_0=0$). 
A random partition $\Pi$ of the set ${\mathbb N}$ into disjoint nonempty blocks is defined by 
assigning two generic integers $m$ and $n$ to the same block if and only if the 
$m$th and the $n$th terms of $(U_j)\giv(H_k)$ hit the same subinterval.
The same partition $\Pi$ can be defined directly in terms of $(X_n)$, by decomposing
$Q_d$ in disjoint layers $Q\setminus L_{R_1}, L_{R_1}\setminus L_{R_2},\ldots$.
Clearly, $T_1,T_2,\ldots$ are 
the minimal integers in the blocks of $\Pi$, and $N_n$ is the number of blocks represented on
the first $n$ integers.

\par The construction of $(U_j)\giv (H_k)$ does not impose any constraints on  the law of the sequence $(H_k)$, which
can be an arbitrary nonincreasing sequence (the induced $\Pi$ is then the most general partially exchangeable
partition \cite{PTRF}). With this in mind, we shall take for a while a more general approach and
assume (as in \cite{Sieve}) that 
$H_k=W_1\cdots W_k$, $k=1,2,\ldots$ where
$W_1,W_2,\ldots$ are independent copies of a random variable
$W\in [0,1]$ with
finite logarithmic 
moments 
$$\mu={\mathbb E}\,[-\log W],~~\sigma^2={\rm Var}\,[-\log W].$$

\begin{proposition}\label{gCLT}
For $n\to\infty$, the variable $N_n$ is 
asymptotically Gaussian with
moments
$${\mathbb E}\,[N_n]\sim {1\over \mu}\,\log n\,,~~~{\rm Var}\,[N_n]\sim {\sigma^2\over\mu^{3}}\,\log n\,.$$
We also have the strong law
$$N_n \sim {1\over \mu}\,\log n\,~~~{\rm a.s.}$$
\end{proposition}
\begin{proof} Our strategy is to
show that $N_n$ is close to $K_n:=\max\{k: H_k>1/n\}$. By the renewal theorem
\cite{Feller}
$K_n$ 
is asymptotically Gaussian with the mean $\mu^{-1}\log n$ and the variance $\sigma^2 \mu^{-3}\log n$
because $K_n$ is just the 
number of epochs on $[0,\log n]$ of the renewal process with steps $-\log W_j$. 
\par By the construction of $(U_j)\giv(H_k)$, we have a dichotomy: 
$U_n\in \,]H_{k},H_{k-1}]$ implies that either $U_n$ will enter the transformed sequence
or will get replaced by some $H_i\geq H_{k}$. 
Let $U_{n1}<\ldots <U_{nn}$ be the order statistics of $U_1,\ldots,U_n$.
It follows that 
\begin{itemize}
\item[(i)]   if $U_{nj}>H_k$ then $N_n\leq k+j$,
\item[(ii)] if $U_{nk}<H_k$ then  $N_n\geq k$.
\end{itemize}
Let 
$\xi_n$ be the number of uniform order statistics smaller than  $1/n$.
By definition,
 $H_{K_n+1}<1/n<H_{K_n}$, hence
$K_n$ and $\xi_n$ are independent and $\xi_n$ is binomial$(n,1/n)$.
By (i), we have $N_n\leq K_n+\xi_n$ where $\xi_n$ is approximately Poisson$(1)$, which yields
the desired upper bound.

\par  Now 
consider
the threshold $s_n=(\log n)^2/n$ and let $J_n:=\max\{k:\,H_k>s_n\}$. 
By (ii), if the number of order statistics smaller than $s_n$ is at least $J_n$
then $N_n\geq J_n$.
Because $\log n\sim \log n-2\log\log n$ 
the index $J_n$  
is still asymptotically Gaussian with the same moments as 
$K_n$. 
On the other hand,
the number of order statistics smaller than $s_n$ is asymptotically Gaussian 
with moments about $(\log n)^2$.
Hence 
elementary  large deviation bounds 
imply that $N_n\geq J_n$ with probability very close to one.
This yields a suitable lower bound, hence the CLT.
 Along the same lines, the strong law of large numbers follows from $N_n\sim K_n$.
\end{proof}

\par Similar limit theorems have been proved by other methods
for the number of blocks of {\it exchangeable} 
partition in \cite{Sieve}, and for a random (size-biased) path in a quad-tree
\cite{Devroye}.

\par Proposition \ref{P1} follows as an instance of Proposition \ref{gCLT} by computing the logarithmic moments as

$$\mu=\ex\,[-\log W]=-g'(0)=d\,,~~\sigma^2={\rm Var}\,[-\log W]=g''(0)-g'(0)^2=d\,.$$

\section{Poisson-paced records}

The probability $p_{n}$ of a chain record at index $n$ is 
equal to the mean height of the last chain record before $n$.
Asymptotics for these quantities follow most easily by poissonisation.

\par Let $(\tau_n)$ be the increasing sequence of points of a homogeneous Poisson point process (PPP) on ${\mathbb R}_+$,
independent of the marks $(X_n)$. The sequence $((X_n,\tau_n),\,n=1,2,\ldots)$ is then 
the sequence of points of a homogeneous PPP in $Q_d\times{\mathbb R}_+$ in the order of increase  of
the time component, which now assumes values in the continuous range ${\mathbb R}_+$.
Let $\widehat{N}_t$ be the number of chain records  and $B_t$ the height of the last chain record on $[0,t]$, that
is
$$\widehat{N}_t=\max\{k: \tau_{T_k}<t\}\,,~~~B_t=H_{\widehat{N}_t}.$$
Clearly, 
$(B_t)$ is the predictable compensator for $(\widehat{N}_t)$, in particular
$${\mathbb E}\, \left[\int_0^t B_s{\rm d}s\right]={\mathbb E}\,[\widehat{N}_t]\,.$$ 
Proposition \ref{gCLT} translates literally as a CLT for $\widehat{N}_t$ as $t\to\infty$.

\par The process $(B_t)$ is Markov time-homogeneous  with a very simple  type of behaviour. 
Given $B_t=b$ the process remains in state $b$ for some rate-$b$ exponential time and then 
jumps to a new state $bW$, with $W$ a stereotypical copy of $H_1$. 
Immediate from this description is the following self-similarity property:
the  law of $(B_t)$ with initial state $B_0=b$ is the same as the law of the process $(bB_{bt})$ with 
$B_0=1$. This kind of process is well defined for arbitrary initial state $b>0$.
See \cite{BC} for features of this process related to the classical
records and
\cite{BerSurv} for more general self-similar (also called semi-stable) 
processes related to increasing L{\'e}vy processes.
The process $(B_t)$ with $B_0=b$ is naturally associated with 
the chain records defined in terms of a homogeneous PPP in  $b\,Q_d\times {\mathbb R}_+$,
with $b\,Q_d$ being the cube with side $[0,b]$.

\par By the self-similarity of $(B_t)$ the moments 
  $$m_\beta(t):={\mathbb E}\,[B_t^\beta]\,$$  
satisfy a renewal-type equation
$$m_{\beta}'(t)=-m_{\beta}(t)+\ex\,[W^\beta m_\beta(tW)]\,.$$
The series solution to  this equation with  the initial value $m_\beta(0)=1$  is
$$m_\beta(t) = \sum_{k=0}^\infty {(-t)^k\over k!}\prod_{j=0}^{k-1} (1-g(j+\beta))~~~{\rm with}~~g(\lambda)=
{1\over (\lambda+1)^d}\,,
$$
as one can check by direct substitution
(see e.g. \cite{BerGn}).

\par Since $m_1(t)$ is the probability that the first arrival after $t$ is a chain record, we have 
the poissonisation identity
$$m_1(t)=e^{-t}\sum_{n=0}^\infty \,{t^n\over n!}\,p_{n+1}\,,$$
which implies, upon equating  coefficients of the series,
\eq
\label{pn}
p_{n}=\sum_{k=0}^{n-1} {n-1\choose k} (-1)^k\prod_{j=0}^{k-1}(1-g(j+1))\,.
\en
This can be compared with the analogous formulas 
$$\underline{p}_n=g(n-1),~~~~~~~ \overline{p}_n=\sum_{k=0}^{n-1}{n-1\choose k}(-1)^k g(k)$$
for the occurencies of strong and weak records, so it would be nice to have a direct combinatorial argument for
(\ref{pn}).

\par For $d=1$ we obtain from (\ref{pn})  the familiar $p_n=1/n$, and for $d=2$ we obtain (surprisingly simple) 
$p_n=1/(2n)$ (for $n>1$).
For $d>2$ the formulas for $p_n$ do not simplify.

\par Factoring
$$1-g(j+\beta)=\prod_{r=1}^d {j+\beta+1-e^{2\pi{\rm i}r /d}\over j+\beta+1}$$
we see that the series for  $m_\beta$ is a generalised hypergeometric function of the type $_dF_d$. 
Exploiting the 
asymptotic properties of this class of functions,
we  determine the asymptotics
as
\eq \label{inmom}
\lim_{t\to\infty} m_\beta(t) t^{\beta}= {1\over -g'(0)} \prod_{r=1}^{\beta-1} {r\over 1-g(r)}= 
{(\beta !)^{d+1}\over \beta d}\prod_{r=2}^\beta {1\over r^d-1}.
\en
where $\beta=1,2,\ldots$ 
Full asymptotic expansion is obtainable in a similar way, see
\cite{BerGn} for details of the method and references.
The depoissonisation of the $\beta=1$ instance implies, quite expectedly,
$$p_n\sim {1\over d\,n}\,,~~~{\rm as~~}n\to\infty\,.$$

\par The following asymptotics for $B_t$ is also derived from (\ref{inmom}) by application of the method of moments.
\begin{proposition}
The random variable $tB_t$ converges, as $t\to\infty$, in distribution and with all moments
to a random variable $Y$ whose moments are given by 
\eq\label{momen}
{\mathbb E}\,[Y^\beta]=
{(\beta !)^{d+1}\over \beta d}\prod_{r=2}^\beta {1\over r^d-1}\,,~~~\beta=1,2,\ldots
\en
\end{proposition}

\vskip0.5cm
\par The law of $Y$, determined uniquely by the moments (\ref{momen}), may be considered as a kind of extreme-value distribution.
In the case $d=1$ we recover well-known $Y\ed E$ with $E$ standard exponential, and 
for $d=2$ we get
$Y\ed EU$ with $E$ and $U$ independent  exponential and uniform random variables.
In general, there is a series representation 
$$Y\ed E_0W_0+\sum_{k=1}^\infty E_k\prod_{j=0}^k W_j$$
where $E_k$'s  are exponential, $W_j$'s for $j>0$ are as before, 
$W_0$ has density 
\begin{equation}\label{stat}
{\mathbb P}(W_0\in {\rm d}s)={\prob(W\leq s)\over sd}\,\,{\rm d}s\,,~~s\in [0,1]
\end{equation}
(which is density of the stationary distribution for the  stick-breaking with factor $W$) and all variables are independent.
Also, $Y$ may be interpreted as an exponential functional of a stationary compound Poisson  process
with initial state $-\log W_0$ and a generic jump  $-\log W$,
see \cite{Exp}.
In the discrete-time setting, 
the same limit law applies to the height of the last chain record 
before $n$.

\section{Scaling limits}

Let $b>0$ be a scaling parameter which we will send to $\infty$.
In the case of one dimension   the point process $\{R_k\}$ of  record values is a self-similar PPP on ${\mathbb R}_+$
with intensity ${\rm d}x/x\,$ (restricted to $x\in [0,1]$).  
The same limit  appears also for the point process of record times $\{T_k/b\}$.
The bivariate point process $\{(bR_k,\,T_k/b), k=1,2,\ldots\}$ has a joint scaling limit 
which may be 
identified with the set of minimal points (the Pareto boundary) of the homogeneous PPP in ${\mathbb R}_+^2$.
See \cite{Nevzorov, Resnick} for these classical results.

\par These facts can be generalised to chain records in $d>1$ dimensions.
Observe that for the values of chain records
we have the component-wise representation 
$$R_k^{(j)}\ed U_1^{(j)}\cdots U_k^{(j)}\,,~~~j=1,\ldots,d\,;~~k=1,2,\ldots$$
with independent uniform $U_k^{(j)}$'s. Therefore, each marginal process
$\{bR_k^{(j)}, k=1,2,\ldots\}$ converges to the same self-similar PPP on ${\mathbb R}_+$.
The vector point process $\{bR_k\}$ converges, as $b\to\infty$, 
to a degenerate limit in ${\mathbb R}_+^d$ which lives on the union of the coordinate axis
(this follows because any level $c/b$ is surpassed by one 
of the marginal $\{R_k^{(j)}\}$'s considerably before the others).
More interestingly, there is a planar limit for the joint process of heights and record times.

\begin{proposition}\label{P3}
The scaled point process $\{(bH_k,T_k/b), k=1,2,\ldots\}$ has a weak  limit as $b\to\infty$.
The limiting point process $\cal R$ in ${\mathbb R}^2_+$ 
is invariant under hyperbolic shifts $(s,t)\mapsto(b s,t/b)$~(with $b>0$),
and the coordinate projections of $\cal R$ are self-similar point processes.
\end{proposition}
\begin{proof} The existence of the limit follows from the analogous result for Poisson-paced marks, and
in the latter setup the result folows from
\cite[Theorem 1]{BC} which, adapted in our framework, guarantees existence of 
the entrance law from $\infty$ 
for the process $(B_t)$ started at $B_0=b$, as $b\to\infty$. 
The hyperbolic invariance follows from  self-similarity of $(B_t)$.
\end{proof}

\par A more explicit construction of  ${\cal R}$ is the following. 
Let ${\cal H}$ be the multiplicatively stationary (that is, self-similar) multiplicative renewal process 
with a generic factor $W$. We may view
$\cal H$ as an extension to ${\mathbb R}_+$ from $[0,1]$ of the
stick-breaking point process 
$\{W_0, W_0W_1, W_0W_1W_2,\ldots\}$ where $W_k\ed W$ and $W_0$ has the stationary density (\ref{stat}). 
Let $\{\xi_k, ~k\in {\mathbb Z}\}$ be the points of $\cal H$ which may be labelled 
 so that $\xi_0=W_0$ is the maximum point of ${\cal H}\cap[0,1]$, and
$\xi_{-1}>1$.
Assign to each $\xi_k$ an arrival time $\sigma_k:=\sum_{i=-\infty}^k E_i/\xi_i$ where the $E_i$'s 
are independent  standard exponential variables, also independent of $\cal H$.
Then  let ${\cal R}:=\{(\xi_k, \sigma_k), ~k\in {\mathbb Z}\}$.
The hyperbolic invariance of $\cal R$  is obvious from the construction and self-similarity of
$\cal H$.
\par The limit process of heights $\cal H$ is not Poisson,
since the law of $W$ is not beta$(\theta,1)$ (for some $\theta>0$). For a similar reason,
the limit process of record times, which is the time-projection of $\cal R$, is also diferent from a Poisson process.
In the discrete-time setting, the dependence of occurencies of chain records 
follows from our interpretaion of chain records in terms of partition $\Pi$ and
a characterisation of the Ewens partitions in \cite{Nacu}
(where it is shown that the independence would force $W$ to be beta$(\theta,1)$).

\par As noticed by Charles Goldie, 
the component-wise logarithmic transform 
$$(-\log(R_k^{(1)}),\ldots,-\log(R_k^{(d)})),~~k=1,2,\ldots\,$$
sends the chain records in $Q_d$ to the sequence
of sites visited by 
a $d$-dimensional random walk whose components are independent 
one-dimensional random walks with
exponentially distributed increments. 
Equivalently, one can consider the upper chain records from the product  exponential distribution in $d$ dimensions.  
In this regime, 
subject to a suitable normalisation,
the values of chain records concentrate near
 the diagonal of the positive orthant.
\vskip0.5cm

\noindent
\centerline{\bf Acknowledgement} 
\vskip0.2cm
The author is indebted to Charles Goldie for the last remark, 
comments on an earlier draft of the paper and
his most stimulating 
interest.

\end{document}